\documentclass[a4paper,11pt]{article}
\usepackage{amsmath, amsfonts, amscd, amssymb, amsthm, enumerate, mathdots}

\newcommand{\Diag}{\operatorname{Diag}}
\newcommand{\Mat}{\operatorname{M}}
\newcommand{\GL}{\operatorname{GL}}
\newcommand{\Ker}{\operatorname{Ker}}
\newcommand{\im}{\operatorname{Im}}
\newcommand{\rk}{\operatorname{rk}}

\def\K{\mathbb{K}}

\theoremstyle{definition}

\theoremstyle{plain}
\newtheorem{theo}{Theorem}

\newtheorem{cor}[theo]{Corollary}

\newtheorem{step}{Step}

\theoremstyle{plain}

\theoremstyle{remark}

\title{The semigroup generated by the similarity class of a singular matrix}
\author{Cl\'ement de Seguins Pazzis\footnote{Lyc\'ee Priv\'e Sainte-Genevi\`eve, 2, rue
de l'\'Ecole des Postes, 78029 Versailles Cedex, FRANCE.}
\footnote{e-mail: dsp.prof@gmail.com}}

\begin{document}
\maketitle
\begin{abstract}
Let $A$ be a singular matrix of $\Mat_n(\K)$, where $\K$ is an arbitrary field.
Using canonical forms, we give a new proof that the sub-semigroup of $(\Mat_n(\K),\times)$ generated by the similarity class of $A$ is the
set of matrices of $\Mat_n(\K)$ with a rank lesser than or equal to that of $A$.
\end{abstract}

\vskip 2mm
\noindent
\emph{AMS Classification:} Primary 15A30. Secondary 15A23.

\vskip 2mm
\noindent
\emph{Keywords:} semigroup, similarity class, companion matrices, rational canonical form, Jordan canonical form.

\section{Introduction}

Let $\K$ be an arbitrary (commutative) field, and denote by $\Mat_n(\K)$ the algebra of square matrices with $n$ rows
and entries in $\K$, and by $\GL_n(\K)$ its group of invertible elements.

Our starting point is the famous theorem of J.A. Erdos
\cite{Erdos}, which states that every singular matrix of $\Mat_n(\K)$
is a product of idempotent ones (for alternative proofs, see \cite{Araujo2}, \cite{Ballantine}, \cite{Djokovic} and \cite{Kim}; for
an extension to principal ideal domains, see \cite{Bhaskara}).
Erdos's theorem and Howie's set-theoretic counterpart \cite{Howie} have had a significant impact
in inspiring the definition of independence algebras (see \cite{Gould} and \cite{Fountain}).

A trivial consequence of Erdos's theorem is that every singular matrix of $\Mat_n(\K)$
is a product of rank $n-1$ idempotent matrices. Since the rank $n-1$ idempotent matrices form a single similarity class,
a natural extension of the question is to determine the sub-semigroup of $(\Mat_n(\K),\times)$
generated by the similarity class of a given singular matrix $A$. It is straightforward
that this semigroup is included in the semigroup consisting of the matrices $M \in \Mat_n(\K)$ such that $\rk M \leq \rk A$.
Remarkably enough, the converse also holds, as was originally proven independently by Ara\'ujo and Silva \cite{Araujo1}
and Grunenfelder, Omladi\v c, Radjavi and Sourour \cite{Grunenfelder}.
The purpose of this note is to give an efficient and relatively elementary proof of that statement which works regardlessly of the ground field
(unlike Ara\'ujo and Silva's proof) and does not rely on prior results about semigroups generated by the similarity class
of an invertible matrix (unlike Grunenfelder et al's proof).
Let us first restate the main theorem:

\begin{theo}\label{main}
Let $A$ be a rank $p$ singular matrix of $\Mat_n(\K)$. Then the semi-group generated by the similarity class of $A$
is
$$S_p:=\bigl\{M \in \Mat_n(\K) : \; \rk M \leq p\bigr\}.$$
\end{theo}

Notice the straightforward corollary:

\begin{cor}
The sub-semigroups of $\Mat_n(\K)$ which are invariant under conjugation and contain only singular matrices are the
$S_k$'s for $k \in \{0,\dots,n-1\}$.
\end{cor}

Before explaining our proof, we need to recall a few standard notations:
given a monic polynomial $P=t^p-\underset{k=0}{\overset{p-1}{\sum}} a_k t^k$ (with indeterminate $t$),
we denote its companion matrix by
$$C(P):=\begin{bmatrix}
0 &   & & 0 & a_0 \\
1 & 0 & &   & a_1 \\
0 & \ddots & \ddots & & \vdots \\
\vdots & & & 0 & a_{p-2} \\
0 & & &  1 & a_{p-1}
\end{bmatrix} \in \Mat_p(\K).$$
We also introduce the (nilpotent) Jordan matrix:
$$J_p:=C(t^p)^T=\begin{bmatrix}
0 & 1 & 0 & \cdots  & 0 \\
0 & 0 & 1 & &  \vdots \\
  &  & \ddots & \ddots &  0 \\
\vdots  &  &    &  0 & 1 \\
0  & \cdots &     & & 0
\end{bmatrix} \in \Mat_p(\K).$$
Similarity of two matrices $A$ and $B$ will be denoted by $A \sim B$, and
the block-diagonal matrix with diagonal blocks $A_1,\dots,A_d$ will be denoted by $\Diag(A_1,\dots,A_d)$.

\section{Proof of the theorem}

Let $A \in \Mat_n(\K)$ with rank $p<n$ and denote by $S$ the semigroup generated by the similarity class of $A$.
Note that $S$ is invariant under conjugation and $S \subset S_p$. Our goal is to show that $S=S_p$.
There are three steps:

\begin{step}
$S$ contains a rank $p$ matrix $M$ such that $\Ker M \oplus \im M=\K^n$ (i.e.\ for which $0$ is a semi-simple eigenvalue).
\end{step}

\begin{step}
$S$ contains a rank $p$ idempotent matrix.
\end{step}

\begin{step}
$S$ contains every rank $p$ matrix $M$ for which $0$ is a semi-simple eigenvalue.
\end{step}

\begin{step}
$S$ contains every rank $p$ matrix of $\Mat_n(\K)$.
\end{step}

\begin{step}
$S$ contains every matrix $N$ of $\Mat_n(\K)$ such that $\rk N <p$.
\end{step}

\subsection{Proof of Step 1}

Using the Fitting decomposition\footnote{See for example \cite{Benson}.}
and the Jordan canonical form\footnote{See for example \cite{Gantmacher}} for the nilpotent part in this decomposition,
one finds a non-singular matrix $Q \in \GL_q(\K)$ and positive integers $i_1,\dots,i_N$ such that
$$A \sim \Diag(Q,J_{i_1},\dots,J_{i_N}).$$
However $J_k \sim J_k^T$ for every positive integer $k$ (simply conjugate $J_k$ with a well-chosen permutation matrix)
hence $S$ contains the product
$$B=\Diag(Q,J_{i_1},\dots,J_{i_N})\,\Diag(Q,J_{i_1}^T,\dots,J_{i_N}^T)= \Diag(Q^2,I_{i_1-1},0,I_{i_2-1},0\dots,I_{i_N-1},0).$$
Since $Q^2$ is non-singular, one has $\rk B=n-N=\rk A$ and obviously $\im B \oplus \Ker B=\K^n$.

\subsection{Proof of Step 2}

Amongst the rank $p$ matrices $M \in S$ satisfying $\im M\oplus \Ker M=E$,
we choose one for which $\rk(M-I_n)$ is minimal. We claim that $M$ is idempotent:
to prove this, it suffices to show that $\dim \Ker(M-I_n)=p$.
We perform a \emph{reductio ad absurdum} and assume $\dim \Ker(M-I_n)<p$.
We obviously lose no generality assuming that
$$M=\Diag(N,0_{n-p}) \quad \textrm{for some $N \in \GL_p(\K)$.}$$
In the rational canonical form of $N$, one of the elementary polynomials is not $t-1$:
we choose such a polynomial $P(t)$, denote by $d$ its degree, and we may then find a non-singular matrix $Q$ such that
$$N \sim \Diag\bigl(Q\,,\,C(P(t))\bigr).$$
Notice that $P(0)\neq 0$,
hence the matrix $N':=\Diag\bigl(C(P(t)),0\bigr)$ of $\Mat_{d+1}(\K)$ has $t\,P(t)$ as minimal polynomial and is therefore similar to $C(t\,P(t))$.
We thus lose no generality assuming that
$$M=\Diag\bigl(Q\,,\,C(t\,P(t))\,,\,0_{n-p-1}) \quad \text{(where the last diagonal block may be empty).}$$
Now, set
$$K=\begin{bmatrix}
(0) & & 1 \\
 & \iddots & \\
1 & & (0)
\end{bmatrix}=(\delta_{i+j,d+2})_{1 \leq i,j \leq d+1}.$$
A straightforward computation yields:
$$K \,C(t\,P(t))\,K^{-1}\,C(t\,P(t))=\begin{bmatrix}
I_d & C \\
0 & 0
\end{bmatrix}$$
for some column-matrix $C \in \Mat_{d+1,1}(\K)$, hence
$$K \,C(t\,P(t))\,K^{-1}\,C(t\,P(t)) \sim \begin{bmatrix}
I_d & 0 \\
0 & 0
\end{bmatrix}.$$
It follows that $S$ contains the product
$$M':=\Diag\bigl(Q\,,\,K C(t\,P(t))K^{-1}\,,\,0_{n-p-1}\bigr)\times
\Diag\bigl(Q\,,\,C(t\,P(t))\,,\,0_{n-p-1}\bigr)
\sim \Diag\bigl(Q^2\,,\,I_d\,,\,0_{n-p}\bigr),$$
for which $0$ is a semi-simple eigenvalue.
The assumptions on $P(t)$ show that $\dim \Ker(C(P(t))-I_d)<d$, whereas $\rk (Q^2-I_n)\leq \rk(Q-I_n)$: 
therefore $\rk(M'-I_n)<\rk(M-I_n)$, which contradicts the definition of $M$ and completes Step 2.

\subsection{Proof of Step 3}

Let $B$ be a rank $p$ matrix of $\Mat_n(\K)$ for which $0$ is a semi-simple eigenvalue.
Hence
$$B \sim \Diag(Q,0_{n-p}) \quad \text{for some $Q \in \GL_p(\K)$.}$$
By Erdos's theorem \cite{Erdos}, the matrix $B':=\Diag(Q,0) \in \Mat_{p+1}(\K)$
may be decomposed as $B'=P_1\cdots P_N$ where the $P_k$'s are rank $p$ idempotent matrices.
Then, for each $k \in \{1,\dots,N\}$, the matrix $Q_k:=\Diag(P_k,0_{n-p-1})$ is a rank $p$ idempotent,
and it therefore belongs to $S$ by the previous step. However
$$Q_1\cdots Q_N=\Diag(Q,0_{n-p}) \sim B$$
hence $B \in S$, which completes the present step.

\subsection{Proof of Step 4}

Let $B \in \Mat_n(\K)$ be a rank $p$ matrix. Set $d:=\dim (\Ker B \cap \im B)$.
If $d=0$, then we know that $B \in S$. Assume that $d \geq 1$.
Using once more the Fitting decomposition and the Jordan canonical form for the nilpotent part of it, we
find an integer $q$, a matrix $P \in \GL_q(\K)$, and integers $i_1,\dots,i_d$ all greater than $1$ such that
$$B \sim \Diag(P\,,\,J_{i_1},\dots,J_{i_d}\,,\,0_{n-p-d}) \sim \Diag(P\,,\,J_{i_1}^T,\dots,J_{i_d}^T\,,\,0_{n-p-d}).$$
For $k \geq 2$, set
$$C_k:=C(t^k-t) \sim \Diag(C(t^{k-1}-1),0) \quad \text{and} \quad C'_k:=\Diag(I_{k-1},0) \in \Mat_k(\K)$$
and notice that these are rank $k-1$ matrices with $0$ as semi-simple eigenvalue. A straightforward computation shows that
$$\forall k \geq 2, \; J_k^T=C_k\,C'_k,$$
hence
$$\Diag(P\,,\,J_{i_1}^T,\dots,J_{i_d}^T\,,\,0_{n-p-d})=\Diag(P\,,\,C_{i_1},\dots,C_{i_d}\,,\,0_{n-p-d})\times
\Diag(I_q\,,\,C'_{i_1},\dots,C'_{i_d}\,,\,0_{n-p-d})$$
and the matrices $\Diag(P\,,\,C_{i_1},\dots,C_{i_d},0_{n-p-d})$ and
$\Diag(I_q\,,\,C'_{i_1},\dots,C'_{i_d},0_{n-p-d})$ have rank $p$ with $0$ as semi-simple eigenvalue.
It follows from Step 3 that both belong to $S$, and hence $B \in S$.

\subsection{Proof of Step 5}

The line of reasoning here is classical but we reproduce it for the sake of completeness.
Let $B \in \Mat_n(\K)$ with rank $r<p$. Then there are non-singular matrices $Q$ and $Q'$ of $\Mat_n(\K)$
such that
$$B=Q\,\Diag(I_r,0_{n-r})\,Q'.$$
For $k \in \{r+1,\dots,p+1\}$, set $D_k:=\Diag(I_{k-1},0,I_{p+1-k},0_{n-p-1})$
and notice that $D_{r+1}\cdots D_{p+1}=\Diag(I_r,0_{n-r})$,
hence
$$B=(QD_{r+1})\,D_{r+2}\,\cdots\,D_{p}\,(D_{p+1}Q')$$
is a product of rank $p$ matrices and therefore belongs to $S$ by Step 4. This completes Step 5, which finishes our proof of Theorem
\ref{main}.

\section{Suggested problems}

We would like to suggest two problems on the issue of semigroups of matrices generated by a similarity class:
\begin{enumerate}[(1)]
\item Describe the semigroup generated by the similarity class of an invertible matrix
with entries in an infinite field. The case where the determinant of the matrix is a root of unity has been solved in \cite{Grunenfelder}
but the solution of the general case is not known yet.
\item Find a truly elementary proof of Theorem \ref{main}, analogous to the one in \cite{Araujo2}.
\end{enumerate}

\section*{Acknowledgements} The author would like to thank the referee for his many thoughtful suggestions.

\end{document}